\documentclass[12pt]{amsart}

\usepackage{amsmath,amssymb,amsfonts,amsthm,amscd,indentfirst}

\usepackage{indentfirst}
\usepackage{hyperref}

\usepackage{color}

\usepackage{setspace}
\usepackage{amsmath, amssymb,amsthm, amsfonts, color}

\theoremstyle{definition}

\theoremstyle{remark}
\newtheorem{remark}{Remark}

\newcommand{\R}{\mathbb R}

\newcommand{\be}{\begin{equation}}
\newcommand{\ee}{\end{equation}}
\newcommand{\bee}{\begin{equation*}}
\newcommand{\eee}{\end{equation*}}

\def\p{\partial}

\def\la{\langle}
\def\ra{\rangle}

\def\Pi{\displaystyle{\mathbb{II}}}

\def\m{\mathfrak{m}}

\def\bCL{\partial  \text{Cube}_{_L}}

\def\iFLp{F^{(i)}_{ +,   L}}
\def\iFLn{F^{(i)}_{-,  L}}

\begin{document}

\title{Measuring Mass via Coordinate Cubes}

\author{Pengzi Miao}
\address[Pengzi Miao]{Department of Mathematics, University of Miami, Coral Gables, FL 33146, USA}
\email{pengzim@math.miami.edu}

\thanks{The author's research was partially supported by NSF grant DMS-1906423.}


\begin{abstract}
Inspired by a formula of Stern that relates scalar curvature to harmonic functions, 
we evaluate the mass of an asymptotically flat $3$-manifold along faces and edges of a large coordinate cube.
In terms of the mean curvature and dihedral angle, the resulting mass formula 
relates to Gromov's  scalar curvature comparison theory for cubic Riemannian polyhedra.
In terms of the geodesic curvature and turning angle of slicing curves, the formula 
realizes the mass as integration of the angle defect detected by the boundary term in the Gauss-Bonnet theorem.
\end{abstract}


\maketitle

\markboth{Pengzi Miao}{Measuring mass via coordinate cubes}

\section{Motivation and mass formulae} 

In \cite{Stern19}, Stern 
gave an intriguing 
formula relating the scalar curvature of a manifold 
to the level set of its harmonic functions. In its simplest form, Stern's formula 
\cite[equation (14)]{Stern19} 
shows 
\be \label{eq-Stern-1}
 \Delta | \nabla u | =   \frac{1}{2   | \nabla u | }   \left[     | \nabla^2 u |^2  +    | \nabla u |^2 ( R - 2 K_{_\Sigma} )  \right]
\ee
near points where $ \nabla u \ne 0$, here 
 $u$ is a harmonic function on a Riemannian $3$-manifold $(M^3, g)$, $ R$ and $K_{_\Sigma}$ denote
the scalar curvature of $g$  and the Gauss curvature of $\Sigma$, the level set of $u$, respectively. 
Applications of the formula to closed manifolds and to compact manifolds with boundary were given by 
Stern \cite{Stern19}, and Stern and Bray \cite{BrayStern19}.

If the manifold $(M^3, g)$ is asymptotically flat,
by applying Stern's formula, Bray,  Kazaras,  Khuri and Stern \cite{BKKS19}  gave a new elegant proof 
of the $3$-dimensional positive mass theorem, which was originally proved by Schoen and Yau \cite{SchoenYau79}, and  Witten \cite{Witten81}. 
Moreover, the result in \cite{BKKS19} provides an explicit lower bound of the mass of $(M, g)$ via a single harmonic function. 

In the context of asymptotically flat manifolds, 
an observation of Bartnik \cite{Bartnik86} was
\be \label{eq-Bartnik-formula}
 \sum_{i=1}^3 \int_{S_\infty}    \frac12 \frac{\p }{\p \nu} | \nabla y^i |^2 \, d \sigma 
  = 16 \pi \m (g), 
\ee
where $\m(g)$ is the mass of $(M,g)$,  $\{ y^i \}$ are harmonic coordinates near infinity,
and $ \int_{S_\infty}$ denotes the limit of integration along a sequence of suitable surfaces  tending  to infinity.
As $ | \nabla y^i | $ approaches $ 1$ sufficiently fast, it can be checked 
\eqref{eq-Bartnik-formula} is equivalent to 
\be \label{eq-Bartnik-formula-1}
 \sum_{i=1}^3  \int_{S_\infty}  \frac{\p }{\p \nu} | \nabla y^i | \, d \sigma = 16 \pi \m (g).
\ee

In view of \eqref{eq-Stern-1} and \eqref{eq-Bartnik-formula-1}, it may be desirable to 
have a formula that computes $\m (g)$
solely in terms of geometric data of the  level sets of $y^i$ near infinity.
In this note, we derive some  formulae of this nature. 
As the level sets of $y^i$ are simply coordinate planes, 
we are thus prompted to  compute $\m(g)$ on the boundary of large coordinate cubes.

A Riemannian $3$-manifold $(M, g)$ is called asymptotically flat with a metric falloff rate $ \tau $ if 
there exists a coordinate chart $ \{ x^i \}$, outside a compact set, 
in which the metric coefficients  satisfies  
\be \label{eq-metric-falloff-c}
g_{ij} = \delta_{ij} + O ( |x |^{-\tau} ), \ 
 \p g_{ij} = O ( |x|^{-\tau -1}) ,  \ \ \p \p g_{ij} = O (|x|^{-\tau -2} ).
\ee
The scalar curvature $R$ of $g$ is assumed  to be integrable so that $\m(g)$ is defined.

 \vspace{.3cm}

\noindent {\bf Geometric Mass Formula.}
{\it Let $ (M^3, g)$ be an asymptotically flat $3$-manifold 
with metric falloff rate $ \tau > \frac12$. Given any large  constant $ L > 0 $, 
let $ \bCL$ denote the boundary of a coordinate cube with side length $2L $ centered at the coordinate origin.
Let $H$ be the mean curvature of the face of $\bCL$ with respect to the outward unit normal $\nu$ in $(M, g)$.
Let $ \mathcal{E}_L$ be the set of all edges of  $ \bCL$. Along each edge in  $ \mathcal{E}_L$, 
let $ \theta$ be the angle between $\nu$ on the two adjacent faces. 
Then, as $L \to \infty$, 
\be \label{eq-Gromov-mass}
\m (g) =  -\frac{1}{8\pi} \int_{\bCL} H \, d \sigma + \frac{1}{8\pi} \int_{\mathcal{E}_L} \left( \frac{\pi}{2} - \theta \right) \, d s + O ( L^{1-2 \tau} ) .
\ee
Here $ d \sigma$ and $ d s$ are the area and the length measure with respect to $g$, respectively. 
Moreover, in terms of the curve $C^{(k)}_t$ which is the  intersection of $\bCL$ and
the coordinate plane $x^k = t$, 
\be \label{eq-Gauss-Bonnet-mass}
\begin{split}
 \m(g)  
=  \frac{1}{8\pi}  \sum_{k=1}^3   \int_{-L}^L 
\left(   2 \pi - \beta^{ \la k \ra }_t 
 -     \int_{ C^{ (k)}_t   } \kappa^{(k)}  \, d s  \right)  \, dt   + O (L^{1 - 2\tau} )   .
\end{split}
\ee
Here $\kappa^{(k)} $ is the geodesic curvature  of $C^{(k)}_t $
and $\beta^{\la k \ra }_t$  is the sum of the turning angle of $C^{(k)}_t$  at its four vertices.
}

\vspace{.2cm}

We give a few remarks regarding these formulae. 

\begin{remark} 
Though our discussion is motivated by \eqref{eq-Stern-1} and  \eqref{eq-Bartnik-formula-1}, the above formulae 
do not assume $\{ x^i \}$ to be harmonic. 
\end{remark}

\begin{remark}
In terms of the dihedral angle $\alpha$ between the two adjacent faces at an edge, 
\eqref{eq-Gromov-mass} is equivalent to 
\be \label{eq-Gromov-mass-1}
\m (g) =  -\frac{1}{8\pi} \int_{\bCL} H \, d \sigma 
+ \frac{1}{8\pi} \int_{\mathcal{E}_L} \left( \alpha - \frac{\pi}{2} \right) \, d s + O ( L^{1-2 \tau} ) .
\ee
In \cite{Gromov14}, Gromov proposed and outlined the proof of a
scalar curvature comparison theorem for polyhedra -- 
let $(D^3, g)$ be a cube-type Riemannian polyhedron with faces $F_j$, let $\alpha_{ij}$ 
be the dihedral angle  between two adjacent faces $F_i$ and $F_j$,  
then the following can not  simultaneously hold:
\begin{itemize}
\item the scalar curvature $R$ of $(D, g)$ is nonnegative; 
\item the mean curvature $H $ of all faces $F_j$ is nonnegative; and
\item the dihedral angle function $\alpha_{ij} < \frac{\pi}{2} $ for all $i$ and $j$.
\end{itemize}
In \cite{Li17}, Li established the corresponding rigidity case under the assumption $\alpha_{ij} \le \frac{\pi}{2} $.
(Further investigation of Gromov's scalar curvature polyhedral comparison theory and edge metrics was given by Li and Mantoulidis \cite{LiMantoulidis18}.)
Now suppose $(M^3, g)$ is a complete, asymptotically flat manifold with nonnegative scalar curvature.
It follows from the positive mass theorem and formula \eqref{eq-Gromov-mass-1} that
\be \label{eq-improve-Gromov}
-\frac{1}{8\pi} \int_{\bCL} H \, d \sigma + \frac{1}{8\pi} \int_{\mathcal{E}_L} \left( \alpha - \frac{\pi}{2} \right) \, d s \ge 0
\ee
for large $L$. These large cubes in $(M, g)$ provide examples for which Gromov's above 
pointwise assumptions on $H$ and $\alpha_{ij}$ may be promoted to an integral inequality. 
\end{remark}

\begin{remark}
Heuristically, if  $\bCL$ could be isometrically embedded in $ \R^3$ as the boundary 
of a standard cube, the right side of \eqref{eq-Gromov-mass} would represent the corresponding 
Brown-York mass of $\bCL$. In this context, formula \eqref{eq-Gromov-mass} resembles the 
convergence of Brown-York mass of  large coordinate spheres to $\m(g)$ (see \cite{FST07}).  
\end{remark}

\begin{remark}
In \eqref{eq-Gauss-Bonnet-mass}, the quantity $ 2 \pi -  \beta^{ \la k \ra }_t  - \int_{ C^{ (k)}_t   } \kappa^{(k)}  \, d s $
measures the angle defect of the large portion of the coordinate plane $\{ x^k = t \}$ inside the cube,
and \eqref{eq-Gauss-Bonnet-mass} shows the mass of $(M^3, g)$  equals 
suitable integration of this angle defect associated to all coordinate planes.
(In the setting of asymptotically conical surfaces, the angle defect can be interpreted as the $2$-d ``mass" of 
those surfaces,  for instance see \cite{CarlottoLellis17}.)
\end{remark}

\begin{remark}
Formulae \eqref{eq-Gauss-Bonnet-mass} is different from the mass formula 
of Bray-Kazaras-Khuri-Stern \cite[equation(6.27)]{BKKS19}.
We will examine this difference  in Section \ref{sec-BKKS-mass}.

If $\{ x^i \}$ are harmonic coordinates,  then upon integration and applying the Gauss-Bonnet theorem,\
Stern's formula  \eqref{eq-Stern-1}, \eqref{eq-Bartnik-formula-1}  and \eqref{eq-Gauss-Bonnet-mass} imply a lower bound of $\m(g)$ 
in the same manner as in \cite{BKKS19}.
For instance,  suppose $M $ has no boundary, consider 
$$U = (u^1, u^2, u^3): (M^3, g) \rightarrow (\R^3, g_0) $$ 
to be a harmonic map, which is a diffeomorphism near infinity 
such that $ g - U^*(g_0)$ satisfies the metric decay condition \eqref{eq-metric-falloff-c}, here $g_0$ is the Euclidean metric. (By the construction 
of harmonic coordinates, for instance in \cite{Bartnik86, Chrusciel86}, this map $U$ always exists.)
Suppose  the regular level set $\Sigma_t^{(i)}$ of all the $u^i$ is connected so that 
$\chi (\Sigma_t^{(i)}) \le 1 $, for instance if $ M $ is $ \R^3$,  
then  it follows from  \cite[equation (14)]{Stern19}, \eqref{eq-Bartnik-formula-1}, 
\eqref{eq-Gauss-Bonnet-mass} and the Gauss-Bonnet theorem that
\be \label{eq-PMT-refined}
\begin{split}
& \  24 \pi \, \m (g) \\
= & \  \lim_{L \to \infty} \sum_{k=1}^3 \int_{\bCL} \frac{\p}{\p \nu} | \nabla u^k | \, d \sigma  +  
 \lim_{L \to \infty}  \sum_{k=1}^3  \int_{-L}^L  \left[    ( 2 \pi - \beta^{ \la k \ra }_t ) - \int_{ C^{ (k)}_t   } \kappa^{(k)}  \, d s  \right]  \, dt  \\
 \ge & \ \sum_{i=1}^3 \ \int_{M} \frac{1}{2}  \left[ \frac{1}{  | \nabla u^i | }  | \nabla^2 u^i |^2 + R |\nabla u^i | \right] d V.
 \end{split}
\ee
We emphasize that \eqref{eq-PMT-refined} is weaker  than the theorem of Bray-Kazaras-Khuri-Stern
 \cite{BKKS19}, because the bound of $\m(g)$ in \cite{BKKS19}  needs only a single harmonic function. 
\end{remark}

\vspace{.2cm}

\noindent {\em Acknowledgements.}  
I want to thank Hubert Bray for encouraging me in the pursuit of the use of \eqref{eq-Bartnik-formula-1}.

\section{Calculation on the cubic boundary}

We will verify \eqref{eq-Gromov-mass} and \eqref{eq-Gauss-Bonnet-mass} by elementary calculation. 
Let $ \{ x^i \}$ be a coordinate chart  of  $(M^3, g)$, outside a compact set,  
in which \eqref{eq-metric-falloff-c} holds.
Given a large  constant $ L > 0 $, 
let $ \bCL$ be the boundary of the coordinate cube with side length $2L $ centered at the coordinate origin.
More precisely, for each $i \in \{ 1, 2, 3 \}$ and $\{j, k \} = \{ 1, 2, 3 \} \setminus \{ i \}$, define the faces 
\bee
\iFLp=  \{ (x^1, x^2, x^3 ) \ | \ x^i = L, | x^j | \le L , | x^k | \le L   \} ,
\eee
\bee
\iFLn=  \{ (x^1, x^2, x^3 ) \ | \ x^i = - L, | x^j | \le L , | x^k | \le L \}.
\eee
Then
\bee
\bCL = \cup_{i=1}^3 \left( F^{i}_{_{+,L} } \cup F^{i}_{_{-,L} } \right).
\eee
For any $ i \neq j $, define the edges 
\bee
E^{(i j )}_{ +, +, L} =
F^{(i)}_{ +, L} \cap F^{(j)}_{ +,  L} , \ 
E^{(i j )}_{ +, -, L} =
F^{(i)}_{ +, L} \cap F^{(j)}_{ -,  L} ,
\eee
and
\bee
E^{(i j )}_{ -, +, L} =
F^{(i)}_{ -, L} \cap F^{(j)}_{ +,  L} , \ 
E^{(i j )}_{ -, -, L} =
F^{(i)}_{ -, L} \cap F^{(j)}_{ -,  L} .
\eee
Let $ \nu $ denote the outward unit $g$-normal to $ \bCL$. Then
\be
\nu = 
\left\{ 
\begin{array}{cc}
\frac{\nabla x^i}{ | \nabla x^i |}  & \ \mathrm{on} \ F^{(i)}_{ +, L}  \\
- \frac{\nabla x^i}{ | \nabla x^i |}  & \ \mathrm{on} \ F^{(i)}_{ -, L}  .
\end{array}
\right.
\ee
Along the edge $E^{(i j )}_{ +, +, L}$, let 
$ \theta^{(i j )}_{ +, +, L} $
be the angle between $\nu$ on the two adjacent faces. 
Then
\be \label{eq-angle-1}
\begin{split}
\cos \theta^{(i j )}_{ +, +, L} = & \  \la \frac{\nabla x^i}{ | \nabla x^i |}  , \frac{\nabla x^j}{ | \nabla x^j |} \ra \\
= & \   | \nabla x^i |^{-1} | \nabla x^j |^{-1} g^{ij} \\
= & \ ( 1 + O ( L^{-\tau} ) ) ( - g_{ij} + O ( L^{- 2\tau} ) ) \\
= & \ - g_{ij}  + O ( L^{- 2 \tau} ),
\end{split}
\ee
where we used the fact 
$ g^{ij} = - g_{ij} + O ( L^{-2 \tau }), \ \text{if} \ i \ne j $.
Similarly, define the angle
$ \theta^{(i j )}_{ +, -, L}$,
$ \theta^{(i j )}_{ -, +, L}$, $ \theta^{(i j )}_{ -, -, L} $
along  the edges $ E^{(i j )}_{ +, -, L}$,
$ E^{(i j )}_{ -, +, L}$, $ E^{(i j )}_{ -, -, L} $, respectively, and we have
\be \label{eq-angle-2}
\begin{split}
\cos \theta^{(i j )}_{ -, +, L} = & \  g_{ij} + O ( L^{- 2 \tau}) \\
\cos \theta^{(i j )}_{ +, -, L} = & \ g_{ij} + O ( L^{- 2 \tau}) \\
\cos \theta^{(i j )}_{ -, -, L} =  & \ - g_{ij} + O ( L^{- 2 \tau}) .
\end{split}
\ee

We are also interested in the intersection between $ \bCL$ and coordinate planes. 
Given any $  t \in [-L, L]$, let 
$ P^{(k)}_t$ denote the coordinate $2$-plane  $ x^k = t $. Let 
$$ C^{(k)}_t = \bCL \cap P^{(k)}_t $$ 
be the ``square" like curve, consisting  of four coordinate curves on the faces  $ F^{(i)}_{\pm, L }$, $ i \ne k$.
Along $C^{(k)}_t$, let $ \kappa^{(k)}$ denote the $g$-geodesic curvature of $ C^{(k)}_t$ 
in $ P^{(k)}_t $ with respect to the outward $g$-unit normal $ \bar \nu$. 

Along  $C^{(k)}_t \cap F^{(i)}_{+, L}$,  $ \bar \nu = \p_{x^i} + O (L^{-\tau})$. Let $ j \in \{ 1, 2, 3 \} \setminus \{k , i \}$, then
\be \label{eq-kappa-p}
\begin{split}
\kappa^{(k)} = & \  - \frac{1}{ g_{jj} }  \la \nabla_{\p_{x^j} } \p_{x^j} , \bar \nu  \ra \\
= & \  - \la \nabla_{\p_{x^j} } \p_{x^j}, \p_{x^i} \ra  + O ( L^{ - 2 \tau - 1} )    \\
= & \ - \Gamma^i_{jj}  + O ( L^{ - 2 \tau - 1} )    \\
= & \  \frac12 g_{jj,i} - g_{ij,j}   + O ( L^{ - 2 \tau - 1} )  .
\end{split}
\ee
Similarly, along  $C^{(k)}_t \cap F^{(i)}_{-, L}$, $ \bar \nu =  - \p_{x^i}  + O (L^{-\tau}) $ and
\be \label{eq-kappa-pp}
\begin{split}
\kappa^{(k)} = & \  - \frac{1}{ g_{jj} }  \la \nabla_{\p_{x^j} } \p_{x^j} , - \p_{x^i}  \ra  + O (L^{- 2 \tau -1} )  \\
= & \  - \left( \frac12 g_{jj,i} - g_{ij,j} \right)  + O ( L^{ - 2 \tau - 1} )  .
\end{split}
\ee

On  $\bCL$,  let $H$ be the mean curvature of its faces in $(M, g)$ with respect to $ \nu$. 
Then, on $ F^{(i)}_{ +, L }$, 
\be \label{eq-H-kappa-sum}
\begin{split}
H = & \ - \sum_{j \ne i, k \ne i } g^{jk} \la \nabla_{\p_{x^j} } \p_{x^k} , \nu \ra \\
= & \ - \sum_{ j \neq i}   \la \nabla_{\p_{x^j} } \p_{x^j} , \p_{x^i} \ra  + O ( L^{ - 2 \tau - 1} )  \\
= & \ \sum_{k \neq i } \kappa^{(k)} + O ( L^{ - 2 \tau - 1} ) .
\end{split}
\ee
Similarly, \eqref{eq-H-kappa-sum} holds on $ F^{(i)}_{ -, L }$ too.

Finally, we measure the turning angle of $C^{(k)}_t $ at each of its vertices. 
At the vertex $ C^{(k)}_t  \cap E^{(ij)}_{+,+, L}$, let $ \beta^{(ij)}_{+,+, L}$ denote the turning angle of $ C^{(k)}_t  $, 
i.e. the angle between $ \p_{x^j}$ and $ - \p_{x^i}$, then 
\be \label{eq-t-angle-1}
\cos \beta^{(ij)}_{+,+, L} = - g_{ji} .
\ee
Similarly, if  $\beta^{(ji)}_{+,-, L}$,  $\beta^{(ij)}_{-,-, L}$, $\beta^{(ji)}_{-,+, L}$ denote 
 the turning angle of $ C^{(k)}_t  $ at  vertices in $E^{(ji)}_{+,-, L}$, $ E^{(ij)}_{-,-, L}$, $ E^{(ji)}_{-,+, L}$, respectively,
 then
\be \label{eq-t-angle-2}
\cos \beta^{(ji)}_{+,-, L} = g_{ij} , \ \cos \beta^{(ij)}_{-,-, L} = - g_{ji} , \ \cos \beta^{(ji)}_{-,+, L} = g_{ij} .
\ee
We define $ \beta^{\la k \ra }_t $ to be  the sum of  the four turning angles of $C^{(k)}_t $ at its vertices.
Then 
\be \label{eq-total-t-angle}
\beta^{\la k \ra }_t = \frac12 \sum_{ \{ i, j \} = \{ 1, 2, 3 \} \setminus \{ k \} } 
\left( \beta^{(ij)}_{+,+, L}  +  \beta^{(ji)}_{+,-, L} +  \beta^{(ij)}_{-,-, L} +  \beta^{(ji)}_{-,+, L} 
\right).
\ee
The factor $\frac12$ here is because of the symmetry 
$ \beta^{(ij)}_{\mu , \lambda, L}  = \beta^{(ji )}_{\lambda, \mu, L}  $
for any indices $i \ne j $ and any sign symbols $\mu, \lambda \in \{ +, - \}$.
(Similarly, $ \theta^{(ij)}_{\mu , \lambda, L}  = \theta^{(ji )}_{\lambda, \mu, L}  $.)

We now turn to the mass $\m(g)$ of $(M^3, g)$. By \cite[Proposition 4.1]{Bartnik86},  $\m(g)$ can be computed by
\be  \label{eq-mass-cubic-int}
\m (g) = \lim_{ L \to \infty} \frac{1}{16 \pi} \int_{\bCL} \sum_{j, k} ( g_{jk,j} -  g_{jj,k} ) \nu^k \, d \sigma .
\ee
Since $ \nu = \p_{x^i} + O ( L^{- \tau} )$ on $ F^{i}_{+, L}$ and $\nu = - \p_{x^i} + O ( L^{- \tau} )$ on $ F^{i}_{-, L}$ , \eqref{eq-mass-cubic-int} simplifies to  
\be \label{eq-mass-cubic-int-1}
\begin{split}
& \ 16 \pi \, \m (g) \\
= & \ \lim_{ L \to \infty}  \sum_{i}  \left(  \int_{F^{i}_{+, L} } \sum_{j \ne i } ( g_{ji,j} -  g_{jj,i} )  \, d \sigma 
 - \int_{F^{i}_{-, L} } \sum_{j \ne i } ( g_{ji,j} -  g_{jj,i} )  \, d \sigma \right) .
 \end{split}
\ee
On $F^{(i)}_{+, L}$, by  \eqref{eq-kappa-p} and \eqref{eq-H-kappa-sum}, 
\be \label{eq-mass-H-1}
\begin{split}
 & \  \int_{F^{(i)}_{ + ,  L}  } \sum_{ j \ne i  } ( g_{ji,j} -  g_{jj,i} ) \, d \sigma  \\
= & \   \int_{F^{(i)}_{ + ,  L}}  \sum_{ j \ne i } \left( - g_{ji,j} \right) \, d \sigma  
-   2 \int_{F^{(i)}_{+, L} } \sum_{k \ne i}   \kappa^{(k)}  \, d \sigma + O ( L^{1-2\tau})  \\
= & \  \int_{F^{(i)}_{ + ,  L}}  \sum_{ j \ne i } ( - g_{ji,j} ) -  2  \int_{F^{(i)}_{ + ,  L}}  H \, d \sigma + O ( L^{1-2\tau})  .
\end{split}
\ee
On each face and edge, let $ d \sigma_{0}$, $ d s_{0}$ denote the area and length measure with respect to the background Euclidean metric $g_0$. 
Then
\be \label{eq-int-gjij-2}
\begin{split}
& \ \int_{F^{(i)}_{ +, L}  }  \sum_{ j \ne i }  g_{ji,j}   \, d \sigma 
=  \int_{F^{(i)}_{ +, L}  }  \sum_{ j \ne i }  g_{ji,j}   \, d \sigma_{0}  + O ( L^{1-2 \tau} ) \\
= & \  \sum_{ j \ne i }  \left[  \int_{E^{(ij)}_{+, +, L} } g_{ji} \, d s_{0} + \int_{E^{(ij)}_{+, -, L} } (-  g_{ji} ) \, d s_{0} \right] + O ( L^{1-2 \tau} ) \\
= & \  \sum_{ j \ne i }  \left[ \int_{E^{(ij)}_{+, +, L} } g_{ji} \, d s + \int_{E^{(ij)}_{+, -, L} } ( - g_{ji} ) \, d s \right] + O ( L^{1-2 \tau} )  \\
= & \  \sum_{ j \ne i }  \left[  - \int_{E^{(ij)}_{+, +, L} }  \cos \theta^{(ij)}_{+, +, L}   \, d s - \int_{E^{(ij)}_{+, -, L} }  \cos \theta^{(ij)}_{+, -, L}    \, d s \right] + O ( L^{1-2 \tau} ), \\
= & \ (-1)  \sum_{ j \ne i }  \left[ \int_{E^{(ij)}_{+, +, L} }  \left(   \frac{\pi}{2} -   \theta^{(ij)}_{+, +, L}  \right)  \, d s + \int_{E^{(ij)}_{+, -, L} }  
 \left(   \frac{\pi}{2}  - \theta^{(ij)}_{+, -, L}  \right)    \, d s  \right] + O ( L^{1-2 \tau} ) ,
\end{split}
\ee
where we have used \eqref{eq-angle-1} and \eqref{eq-angle-2}.
It follows from  \eqref{eq-mass-H-1} and \eqref{eq-int-gjij-2} that
\be \label{eq-mass-fi}
\begin{split}
& \  \int_{F^{(i)}_{ + ,  L}  } \sum_{j \ne i} ( g_{ji,j} -  g_{jj,i} )  \, d \sigma \\
= & \  \sum_{ j \ne i }  \left[ \int_{E^{(ij)}_{+, +, L} }  \left(   \frac{\pi}{2} -   \theta^{(ij)}_{+, +, L}  \right)  \, d s + \int_{E^{(ij)}_{+, -, L} }  
 \left(   \frac{\pi}{2}  - \theta^{(ij)}_{+, -, L}  \right)    \, d s  \right] \\
 & \ -  2  \int_{F^{(i)}_{ + ,  L}}  H \, d \sigma  + O ( L^{1-2 \tau} )  .
\end{split}
\ee
Similarly,  on $ F^{(i)}_{ - ,  L}$, 
\be \label{eq-mass-fi-n}
\begin{split}
& \  \int_{F^{(i)}_{ - ,  L}  } \sum_{j \ne i} ( g_{ji,j} -  g_{jj,i} ) (-1) \, d \sigma \\
= & \  \sum_{ j \ne i }  \left[ \int_{E^{(ij)}_{-, +, L} }  \left(   \frac{\pi}{2} -   \theta^{(ij)}_{-, +, L}  \right)  \, d s + \int_{E^{(ij)}_{-, -, L} }  
 \left(   \frac{\pi}{2}  - \theta^{(ij)}_{-, -, L}  \right)    \, d s  \right] \\
 & \ -  2  \int_{F^{(i)}_{ - ,  L}}  H \, d \sigma  + O ( L^{1-2 \tau} )  .
\end{split}
\ee
By \eqref{eq-mass-cubic-int-1}, \eqref{eq-mass-fi} and \eqref{eq-mass-fi-n}, we have
\be \label{eq-mass-H-angle-1}
\begin{split}
& \ 16 \pi \m(g)  \\
= & \  \sum_{  j \ne i  }  \left[ \int_{E^{(ij)}_{+, +, L} }  \left(   \frac{\pi}{2} -   \theta^{(ij)}_{+, +, L}  \right)  \, d s + \int_{E^{(ij)}_{+, -, L} }  
 \left(   \frac{\pi}{2}  - \theta^{(ij)}_{+, -, L}  \right)    \, d s \right. \\
& \ \left. + \int_{E^{(ij)}_{-, +, L} }  \left(   \frac{\pi}{2} -   \theta^{(ij)}_{-, +, L}  \right)  \, d s + \int_{E^{(ij)}_{-, -, L} }  
 \left(   \frac{\pi}{2}  - \theta^{(ij)}_{-, -, L}  \right)    \, d s
   \right] \\
 & \ -  2  \int_{\bCL}  H \, d \sigma  + O ( L^{1-2 \tau} )  ,
\end{split}
\ee
which verifies \eqref{eq-Gromov-mass}. Note that each edge of $ \bCL$ is 
counted twice in \eqref{eq-mass-H-angle-1}.

Next, we write $\m(g)$ in terms of the geodesic curvature and turning angles of $C^{(k)}_t$.
By \eqref{eq-angle-1}, \eqref{eq-angle-2} and \eqref{eq-t-angle-1}, \eqref{eq-t-angle-2}, 
we have
\be
 \theta^{(ij)}_{\mu, \lambda, L} =  \beta^{(ij)}_{\mu, \lambda, L} + O ( L^{-2\tau}),
\ee
for any indices $i \ne j $ and any sign symbols $\mu, \lambda \in \{ +, - \}$.
Thus, 
\be \label{eq-two-set-angles}
\begin{split}
 & \  \sum_{  j \ne i  }  \left[ \int_{E^{(ij)}_{+, +, L} }  \left(   \frac{\pi}{2} -   \theta^{(ij)}_{+, +, L}  \right)  \, d s + \int_{E^{(ij)}_{+, -, L} }  
 \left(   \frac{\pi}{2}  - \theta^{(ij)}_{+, -, L}  \right)    \, d s \right. \\
& \ \left. + \int_{E^{(ij)}_{-, +, L} }  \left(   \frac{\pi}{2} -   \theta^{(ij)}_{-, +, L}  \right)  \, d s + \int_{E^{(ij)}_{-, -, L} }  
 \left(   \frac{\pi}{2}  - \theta^{(ij)}_{-, -, L}  \right)    \, d s \right] \\
 = & \   \sum_{  j \ne i  }  \left[ \int_{E^{(ij)}_{+, +, L} }  \left(   \frac{\pi}{2} -   \beta^{(ij)}_{+, +, L}  \right)  \, d s_0 + \int_{E^{(ij)}_{+, -, L} }  
 \left(   \frac{\pi}{2}  - \beta^{(ij)}_{+, -, L}  \right)    \, d s_0 \right. \\
& \ \left. + \int_{E^{(ij)}_{-, +, L} }  \left(   \frac{\pi}{2} - \beta^{(ij)}_{-, +, L}  \right)  \, d s_0 + \int_{E^{(ij)}_{-, -, L} }  
 \left(   \frac{\pi}{2}  - \beta^{(ij)}_{-, -, L}  \right)    \, d s_0 \right]  + O ( L^{1-2\tau}) \\
= & \   \sum_k   \int_{-L}^L 
2 \left( 2 \pi - \beta^{ \la k \ra }_t \right)  \, d t   + O ( L^{1-2\tau}) .
\end{split}
\ee

By \eqref{eq-H-kappa-sum}, we also have
\be \label{eq-two-set-curvatures}
\begin{split}
& \ \int_{\bCL} H \, d \sigma 
=  \sum_i  \int_{ F^{(i)}_{+, L} \cup F^{(i)}_{-, L} } H \, d \sigma  \\
= & \ \sum_i  \int_{ F^{(i)}_{+, L} \cup F^{(i)}_{-, L} } \sum_{k \ne i} \kappa^{(k)}  \, d \sigma_0   + O (L^{1- 2 \tau}) \\
= & \ \sum_k   \sum_{i \ne k}   \int_{ F^{(i)}_{+, L} \cup F^{(i)}_{-, L}  } \kappa^{(k)}  \, d \sigma_0  + O (L^{1 - 2\tau} ) \\
= & \ \sum_k   \sum_{i \ne k}  \int_{-L}^L \left(  \int_{ C^{ (k)}_t  \cap \left( F^{(i)}_{+, L}  \cup F^{(i)}_{-, L} \right) } \kappa^{(k)}  \, d s_0 \right) \, dt   + O (L^{1 - 2\tau} )\\
= & \ \sum_k    \int_{-L}^L   \left(  \int_{ C^{ (k)}_t   } \kappa^{(k)}  \, d s_0 \right)  \, dt   + O (L^{1 - 2\tau} ) \\
= & \ \sum_k    \int_{-L}^L   \left(  \int_{ C^{ (k)}_t   } \kappa^{(k)}  \, d s \right)  \, dt   + O (L^{1 - 2\tau} ) .
\end{split}
\ee
Therefore, by  \eqref{eq-mass-H-angle-1}, \eqref{eq-two-set-angles} and \eqref{eq-two-set-curvatures},  we have
\be \label{eq-mass-kappa-angle}
\begin{split}
& \ 16 \pi \m(g)  \\
= & \   2  \sum_k   \int_{-L}^L 
\left(   2 \pi - \beta^{ \la k \ra }_t 
 -     \int_{ C^{ (k)}_t   } \kappa^{(k)}  \, d s  \right)  \, dt   + O (L^{1 - 2\tau} )   ,
\end{split}
\ee
which verifies \eqref{eq-Gauss-Bonnet-mass}.

\section{Relation to the mass formula in \cite{BKKS19}} \label{sec-BKKS-mass}

In formulae \eqref{eq-Gromov-mass} and \eqref{eq-Gauss-Bonnet-mass}, 
the coordinates $\{ x^i \}$ used in defining $\bCL$ and $ C^{(k)}_t$ do not need to be harmonic. 
If  $\{ x^i \}$ are harmonic,  \eqref{eq-Gauss-Bonnet-mass} and \eqref{eq-Bartnik-formula-1} then imply 
\be \label{eq-BM-mass}
\begin{split}
\sum_{k=1}^3  \int_{\bCL} \frac{\p}{\p \nu} | \nabla x^k | d \sigma  +  
 \int_{-L}^L  \left[    ( 2 \pi - \beta^{ \la k \ra }_t ) - \int_{ C^{ (k)}_t   } \kappa^{(k)}  \, d s  \right] dt 
= 24 \pi \m (g) + o(1).
\end{split}
\ee 
This formula  is weaker than that of Bray-Kazaras-Khuri-Stern \cite{BKKS19},   which indicates,
without summing over $k$, each  summand above tends to $8\pi \m (g)$, provided $\Delta x^k =0$.

We  now examine the  summand in \eqref{eq-BM-mass}. 
Let  $k$, $i$, $j$ be fixed indices so that they are distinct from each other. 
Similar to how \eqref{eq-int-gjij-2} is derived, by \eqref{eq-t-angle-1} and  \eqref{eq-t-angle-2}, 
\bee
\begin{split}
2  \int_{-L}^L \left( 2 \pi - \beta^{ \la k \ra }_t \right)  \, d t 
 = & \   - \int_{F^i_{+, L}} g_{ij,j} \, d \sigma_0 +  \int_{F^i_{-, L} } g_{ij, j} \, d \sigma_0 \\
& \  - \int_{F^{j}_{+, L} }  g_{ij,i} \, d \sigma_0   + \int_{F^{j}_{-, L}} g_{ij,i} \, d \sigma_0
   + O ( L^{1 - 2\tau } ) .
\end{split}
\eee
By \eqref{eq-kappa-p} and \eqref{eq-kappa-pp}, 
\bee
\begin{split}
& \ 2 \int_{-L}^L  \left[ -  \int_{ C^{ (k)}_t   } \kappa^{(k)}  \, d s \right]  \, d t \\
= & \   \int_{-L}^L  \left[
\int_{ C^{ (k)}_t  \cap F^i_{+, L}  }  ( 2 g_{ij,j} - g_{jj,i} )  \, d s 
+ \int_{ C^{ (k)}_t  \cap F^i_{-, L}  }  ( - 2 g_{ij,j} + g_{jj,i} )  \, d s \right. \\
& \ \left. + \int_{ C^{ (k)}_t  \cap F^j_{+, L}  }  ( 2 g_{j i,i} - g_{ii,j} )   \, d s 
+ \int_{ C^{ (k)}_t  \cap F^j_{-, L}  }   ( - 2 g_{j i,i} +  g_{ii,j} )    \, d s  \right] \, d t \\
= & \   \int_{  F^i_{+, L}  }  ( 2 g_{ij,j} - g_{jj,i} )  \, d \sigma_0 
+ \int_{ F^i_{-, L}  }  ( - 2 g_{ij,j} + g_{jj,i} )  \, d \sigma_0 \\
& \  + \int_{ F^j_{+, L}  }  ( 2 g_{j i,i} - g_{ii,j} )   \, d \sigma_0
+ \int_{  F^j_{-, L}  }   ( - 2 g_{j i,i} +  g_{ii,j} )    \, d \sigma_0 + O (L^{1 - 2 \tau}) .
\end{split}
\eee
Thus, the boundary term in the Gauss-Bonnet theorem satisfies 
\be \label{eq-ckt-angle}
\begin{split}
& \ 2  \int_{-L}^L 
\left[    ( 2 \pi - \beta^{ \la k \ra }_t ) - \int_{ C^{ (k)}_t   } \kappa^{(k)}  \, d s  \right]  \, dt \\
= & \   \int_{  F^i_{+, L}  }  (  g_{ij,j} - g_{jj,i} )  \, d \sigma_0 
-  \int_{ F^i_{-, L}  }  (  g_{ij,j} - g_{jj,i} )  \, d \sigma_0 \\
& \  + \int_{ F^j_{+, L}  }  (  g_{j i,i} - g_{ii,j} )   \, d \sigma_0
-  \int_{  F^j_{-, L}  }   ( g_{j i,i} -  g_{ii,j} )    \, d \sigma_0 
   + O ( L^{1 - 2\tau } ) .
\end{split}
\ee
(If summing over $ k \in \{ 1, 2, 3\}$, this again gives \eqref{eq-mass-kappa-angle}.)

We next compute $ | \nabla x^k |$ and $ \Delta x^k $.
Since $ | \nabla x^k |^2 = g^{kk}$,
\bee
2 \frac{\p}{\p \nu} | \nabla x^k |
= -  g_{kk, m} \nu^m + O ( L^{-2 \tau -1} ).
\eee
Hence, 
\be
\begin{split}
 2 \int_{\bCL} \frac{\p}{\p \nu} | \nabla x^k | \, d \sigma  & \ 
=  - \int_{F^{k}_{+, L} } g_{kk,k} \, d \sigma_0  + \int_{F^{k}_{-, L} } g_{kk,k} \, d \sigma_0 \\
& \  - \int_{F^{i}_{+, L} } g_{kk,i} \, d \sigma_0  + \int_{F^{i}_{-, L} } g_{kk,i} \, d \sigma_0 \\
& \  - \int_{F^{j}_{+, L} } g_{kk,j} \, d \sigma_0  + \int_{F^{j}_{-, L} } g_{kk,j} \, d \sigma_0 
+ O ( L^{1 - 2 \tau}) .
\end{split}
\ee
The term $g_{kk,k}$ appears in 
\be
\begin{split}
\Delta x^k 
= & \ \sum_{m \ne k } \left( \frac12 g_{mm,k} - g_{km,m} \right) -  \frac12 g_{kk,k} +O ( L^{-2 \tau -1} )
\end{split}.
\ee
Thus, 
\be \label{eq-pu-gdu}
\begin{split}
& \ 2 \int_{\bCL} \frac{\p}{\p \nu} | \nabla x^k | \, d \sigma  - \int_{F^k_{+,L} } 2 \Delta x^k \, d \sigma 
+  \int_{F^k_{-,L} } 2 \Delta x^k  \, d \sigma  \\
= & \  \int_{F^k_{+,L} }   \sum_{m \ne k } \left(  g_{km,m} - g_{mm,k}  \right)  
+ \sum_{m \ne k} g_{km,m}  \, d \sigma_0 \\
& \ -  \int_{F^k_{-,L} } 
 \sum_{m \ne k } \left(  g_{km,m} - g_{mm,k}  \right)  + \sum_{m \ne k} g_{km,m}   \, d \sigma_0 \\
& \  - \int_{F^{i}_{+, L} } g_{kk,i} \, d \sigma_0  + \int_{F^{i}_{-, L} } g_{kk,i} \, d \sigma_0 \\
& \  - \int_{F^{j}_{+, L} } g_{kk,j} \, d \sigma_0  + \int_{F^{j}_{-, L} } g_{kk,j} \, d \sigma_0
+ O ( L^{1 - 2 \tau}) .
\end{split}
\ee
Therefore, adding \eqref{eq-ckt-angle} and \eqref{eq-pu-gdu}, and using \eqref{eq-mass-cubic-int-1}, we have
\be \label{eq-BKKS-mass}
\begin{split}
& \ 2 \int_{\bCL} \frac{\p}{\p \nu} | \nabla x^k | \, d \sigma+  2  \int_{-L}^L 
\left[    ( 2 \pi - \beta^{ \la k \ra }_t ) - \int_{ C^{ (k)}_t   } \kappa^{(k)}  \, d s_g  \right]  \, dt \\
= & \   16 \pi \, \m(g) +  \int_{F^k_{+,L} } 2 \Delta  x^k  \, d \sigma_0 -  \int_{F^k_{-,L} } 2 \Delta  x^k \, d \sigma_0  
  +   O ( L^{1 - 2\tau } )  \\
& \ + \int_{F^{k}_{+, L} }  ( g_{k i,i }  + g_{k j, j}  )    \, d \sigma_0 - \int_{F^k_{-, L} }  ( g_{k i ,i } + g_{k j, j} ) \, d \sigma_0    \\
& \ - \int_{F^i_{+, L} } g_{ik,k} \, d \sigma_0 + \int_{F^i_{-, L} } g_{i k,k} \, d \sigma_0
- \int_{F^j_{+, L} } g_{j k ,k} \, d \sigma_0 + \int_{F^{j}_{-, L} } g_{j k ,k}  \, d \sigma_0  .
\end{split}
\ee
The last two lines in \eqref{eq-BKKS-mass} cancel upon integration by parts. Thus,
\be \label{eq-BKKS-mass-1}
\begin{split}
& \  \int_{\bCL} \frac{\p}{\p \nu} | \nabla x^k | \, d \sigma +    \int_{-L}^L 
\left[    ( 2 \pi - \beta^{ \la k \ra }_t ) - \int_{ C^{ (k)}_t   } \kappa^{(k)}  \, d s  \right]  \, dt \\
= & \   8 \pi \, \m(g)  +   \int_{F^k_{+,L} }  \Delta  x^k  \, d \sigma_0 -  \int_{F^k_{-,L} }  \Delta  x^k \, d \sigma_0  
+  O ( L^{1 - 2\tau } )  ,
\end{split}
\ee
which is the formula in \cite[equation (6.27)]{BKKS19} if $\Delta x^k=0$.

\end{document}